\documentclass[reqno]{amsart}

\newcommand{\Tr}{\text{Tr}} 

\newcommand{\D}{\mathcal{D}}
\newcommand{\bbn}{\protect{\mathbb N}}
\newcommand{\one}{{\bf 1}}
\newcommand{\A}{{\mathcal A}}
\usepackage{amsbsy}

\usepackage{times}
\usepackage[T1]{fontenc}
\usepackage{mathrsfs}
\usepackage{latexsym}
\usepackage[dvips]{graphics}
\usepackage{epsfig}
\usepackage{hyperref, flafter,epsf}
\usepackage{amsmath,amsfonts,amsthm,amssymb,amscd}
\usepackage{color}
\usepackage{fix-cm}
\usepackage{hyperref}
\usepackage{color}
\usepackage{url}
\newcommand{\bburl}[1]{\textcolor{blue}{\url{#1}}}

\usepackage{comment}
\usepackage{multirow}

\newcommand\be{\begin{equation}}
\newcommand\ee{\end{equation}}
\newcommand\bea{\begin{eqnarray}}
\newcommand\eea{\end{eqnarray}}
\newcommand\bi{\begin{itemize}}
\newcommand\ei{\end{itemize}}
\newcommand\ben{\begin{enumerate}}
\newcommand\een{\end{enumerate}}

\newtheorem{thm}{Theorem}[section]
\newtheorem{conj}[thm]{Conjecture}

\newtheorem{lem}[thm]{Lemma}

\theoremstyle{definition}
\newtheorem{exa}[thm]{Example}

\theoremstyle{definition}
\newtheorem{defi}[thm]{Definition}

\newtheorem{rek}[thm]{Remark}

\theoremstyle{definition}
\newtheorem{cla}[thm]{Claim}

\newcommand{\E}{\mathbb{E}}

\newcommand\cycle[2][\,]{%
  \readlist\thecycle{#2}%
  (\foreachitem\i\in\thecycle{\ifnum\icnt=1\else#1\fi\i})%
}

\newcommand{\gl}{\lambda}
\newcommand{\p}{\ensuremath{\mathbb{P}}}

\numberwithin{equation}{section}

\makeatletter
\def\@tocline#1#2#3#4#5#6#7{\relax
  \ifnum #1>\c@tocdepth 
  \else
    \par \addpenalty\@secpenalty\addvspace{#2}%
    \begingroup \hyphenpenalty\@M
    \@ifempty{#4}{%
      \@tempdima\csname r@tocindent\number#1\endcsname\relax
    }{%
      \@tempdima#4\relax
    }%
    \parindent\z@ \leftskip#3\relax \advance\leftskip\@tempdima\relax
    \rightskip\@pnumwidth plus4em \parfillskip-\@pnumwidth
    #5\leavevmode\hskip-\@tempdima
      \ifcase #1
       \or\or \hskip 1em \or \hskip 2em \else \hskip 3em \fi%
      #6\nobreak\relax
    \hfill\hbox to\@pnumwidth{\@tocpagenum{#7}}\par
    \nobreak
    \endgroup
  \fi}
\makeatother

\begin{document}

\title{Distribution of Eigenvalues of Matrix Ensembles arising from Wigner and Palindromic Toeplitz Blocks}

\author[Blackwell]{Keller Blackwell}
\email{kellerb@cs.stanford.edu}
\address{Department of Computer Science, Stanford University, Stanford, CA 94350}

\author[Borade]{Neelima Borade}
\email{nborad2@uic.edu}
\address{Department of Mathematics, University of Illinois at Chicago, Chicago, IL 60607}

\author[Bose]{Arup Bose}
\email{bosearu@gmail.com}
\address{Statistics and Mathematics Unit, Indian Statistical Institute, Kolkata 700108, India}

\author[Devlin VI]{Charles Devlin VI}
\email{chatrick@umich.edu}
\address{Department of Mathematics, University of Michigan, Ann Arbor, MI 48109}

\author[Luntzlara]{Noah Luntzlara}
\email{nluntzla@umich.edu}
\address{Department of Mathematics, University of Michigan, Ann Arbor, MI 48109}

\author[Ma]{Renyuan Ma}
\email{ma.2005@osu.edu}
\address{Department of Mathematics, Ohio State University, Columbus, OH 43210}

\author[Miller]{Steven J. Miller}
\email{sjm1@williams.edu}
\address{Department of Mathematics, Williams College, MA 01267}

\author[Mukherjee]{Soumendu Sundar Mukherjee}
\email{soumendu041@gmail.com}
\address{Interdisciplinary Statistical Research Unit, Indian Statistical Institute, Kolkata 700108, India}

\author[Wang]{Mengxi Wang}
\email{mengxiw@umich.edu}
\address{Department of Mathematics, University of Michigan, Ann Arbor, MI 48109}

\author[Xu]{Wanqiao Xu}
\email{wanqiaox@umich.edu}
\address{Department of Mathematics, University of Michigan, Ann Arbor, MI 48109}



\thanks{This work was supported by NSF grants DMS1561945 and DMS-1659037, the J.C.~Bose National Fellowship, an INSPIRE Faculty fellowship, the University of Michigan, and Williams College; it is a pleasure to thank them for their support. We would also like to thank Shiliang Gao and Dr. Jun Yin for helpful comments, and Zhijie Chen, Jiyoung Kim, and Samuel Murray for providing a counter-example to a conjecture on the contribution of terms in the expansion of formulas for the moments of the disco of two ensembles.}

\maketitle

\begin{abstract} Random Matrix Theory (RMT) has successfully modeled diverse systems, from energy levels of heavy nuclei to zeros of $L$-functions;
this correspondence has allowed RMT to successfully predict many number theoretic behaviors. However there are some operations which to date have no RMT analogue. Our motivation is to find an RMT analogue of Rankin-Selberg convolution, which constructs a new $L$-functions from an input pair. We report one such attempt; while it does not appear to model convolution, it does create new ensembles with properties hybridizing those of its constituents.


For definiteness we concentrate on the ensemble of palindromic real symmetric Toeplitz (PST) matrices and the ensemble of real symmetric matrices, whose limiting spectral measures are the Gaussian and semi-circular distributions, respectively; these were chosen as they are the two extreme cases in terms of moment calculations. For a PST matrix $A$ and a real symmetric matrix $B$, we construct an ensemble of random real symmetric block matrices whose first row is $\lbrace A, B \rbrace$ and whose second row is $\lbrace B, A \rbrace$. By Markov's Method of Moments and the use of free probability, we show this ensemble converges weakly and almost surely to a new, universal distribution with a hybrid of Gaussian and semi-circular behaviors. We extend this construction by considering an iterated concatenation of matrices from an arbitrary pair of random real symmetric sub-ensembles with different limiting spectral measures. We prove that finite iterations converge to new, universal distributions with hybrid behavior, and that infinite iterations converge to the limiting spectral measure of the dominant component matrix.
\end{abstract}

\tableofcontents


\section{Introduction}
\subsection{History}

Random Matrix Theory (RMT) is well-suited to the fundamental problem of studying spacings between observed values arising from large, complex systems such as energy levels of heavy nuclei and vertical spacings of zeros of the Riemann zeta function. Similar to the Central Limit Theorem, the behavior of a typical element is often close to the system average, which frequently can be computed. For example, the entries of Hamiltonians describing the energy levels of heavy nuclei are inextricably complex, but the distribution of the energy levels of these operators are well approximated by the average eigenvalue behavior of the real symmetric ensemble.
This observation is captured by Wigner's Semi-Circular Law \cite{Wig1}, which states that the distribution of normalized eigenvalues of a random real symmetric or complex Hermitian matrix with entries i.i.d.r.v. from a fixed probability distribution with mean $0$ and variance $1$ converges almost surely to the semi-circular density.

The ensemble of $N \times N$ real symmetric matrices has $N(N+1)/2$ independent parameters; a natural question is how placing additional structural constraints, and thereby reducing the degrees of freedom, affects eigenvalue behavior.
Recall an $N\times N$ palindromic symmetric Toeplitz (PST) matrix $A_N$ is of the form
\be\label{eq:defrsptmat}
A_N\ =\ \left(
\begin{array}{ccccccc}
b_0    &  b_1   &  b_2   &  \cdots  &  b_2   &  b_1  &  b_0  \\
b_1    &  b_0   &  b_1   &  \cdots  &  b_3   &  b_2  &  b_1  \\
b_2    &  b_1   &  b_0   &  \cdots  &  b_4   &  b_3  &  b_2  \\
\vdots & \vdots & \vdots &  \ddots  & \vdots & \vdots &  \vdots  \\
b_2    & b_3    & b_4    &  \cdots  & b_0    & b_1   &  b_2   \\
b_1    & b_2    & b_3    &  \cdots  & b_1    & b_0   &  b_1  \\
b_0    & b_1    & b_2    &  \cdots  & b_2    & b_1   &  b_0
\end{array}\right),
\ee
which is a symmetric Toeplitz matrix whose first row is a palindrome.
Bai \cite{Bai} first posed the problem of studying the limiting eigenvalue
distribution associated with random symmetric (non-palindromic) Toeplitz matrices.
Bose-Chatterjee-Gangopadhyay \cite{BCG}, Bryc-Dembo-Jiang \cite{BDJ}, and
Hammond-Miller \cite{HM} independently observed that the limiting
even moments
of random symmetric Toeplitz matrices are
dominated by those of the standard Gaussian.
Subsequent work by Massey, Miller, and Sinsheimer \cite{MMS}
shows that the moments of the PST ensemble
are those of the standard Gaussian, and that the limiting spectral measure converges weakly to the same.

An $N \times N$ real symmetric matrix $B$ has $N(N+1)/2$ degrees of freedom; in contrast, a PST matrix of the same dimensions has only $N/2$ degrees of
freedom. The PST ensemble is then a very thin sub-ensemble of all real symmetric matrices, and the imposed structure leads to new behavior.
Thus by examining sub-ensembles of real symmetric matrices, one has the exciting possibility of seeing new, universal distributions.

We now describe the motivation for the new construction in this paper. The entrance of random matrix theory into number theory came in the '70s in a fortuitous meeting between Hugh Montgomery and Freeman Dyson, yielding the observation that the pair correlation function of Riemann zeta zeros matched that of the eigenvalues of random Hermitian matrices in the Gaussian Unitary Ensemble (see \cite{BFMT-B,FM} for a fuller treatment and history).  Work by Hejhal \cite{Hej} and Rudnick and Sarnak \cite{RS} extended this random matrix connection to $n$-level correlations of zeros of $L$-functions, generalizations of the Riemann zeta function which arise throughout number theory.
The zero densities of $L$-functions can be recast as the study of eigenvalue behavior of random complex Hermitian matrices \cite{FM,Hej,RS}, and Rankin-Selberg convolution allows the creation of a new $L$-function from multiple input $L$-functions. Given families of $L$-functions $\left \lbrace L (s, f_i)_{f_i \in \mathcal{F}_i}\right\rbrace$ with $i \in \lbrace 1, 2, \ldots, I \rbrace)$, the Rankin-Selberg convolution
\begin{equation}
\left \lbrace L ( s, f_1 \otimes \cdots \otimes f_I ) \right \rbrace_{ ( f_1, \ldots, f_I) \in \mathcal{F}_1 \times \cdots \times \mathcal{F}_I }
\end{equation}
gives a new family of $L$ functions\footnote{For details see \cite{IK}.}. Due\~nez and Miller \cite{DM1, DM2} were able to describe the behavior of the zeros of the convolution in terms of the behavior of the constituent families in many situations. As RMT has successfully modeled so many properties of $L$-functions, it is thus natural to ask if there is an RMT analogue of convolutions. Motivated by the confluence
of number theory and random matrix theory,
we consider the eigenvalue behavior of the
ensemble constructed as the ``disco'' concatenation\footnote{The whimsical
naming of the ``disco'' construction arises from the entries of the block
matrix ``ABBA'', a quintessential icon of disco music's heyday.} of PST matrices $A$ and real symmetric matrices $B_1$:

\begin{equation}\label{eq:1-disco-construction}
    \mathcal{D}_1 \left( A, B_1 \right)\ = \ \begin{bmatrix}
    A & B_1 \\
    B_1 & A
    \end{bmatrix}.
\end{equation}

The resulting ensemble of $2N \times 2N$ symmetric block matrices have
only $(N/2) + N(N+1)/2$ degrees of freedom and constitute another thin
subset of all real symmetric matrices that may give rise to new
eigenvalue behavior of interest. We chose the PST and real symmetric ensembles as their limiting distributions
(Gaussian and semi-circular, respectively) demonstrate extreme contrasting behavior: the Gaussian
features a sharp decay rate but unbounded support, while the
semi-circular distribution is strictly bounded within $[-2,2]$. The novel construction from
known ensembles
furthermore poses the question of how the disco ensemble's
limiting eigenvalue distribution may be described in terms
of its constituent distributions.


\begin{figure}[ht]
     \centering
     \includegraphics[scale=0.225]{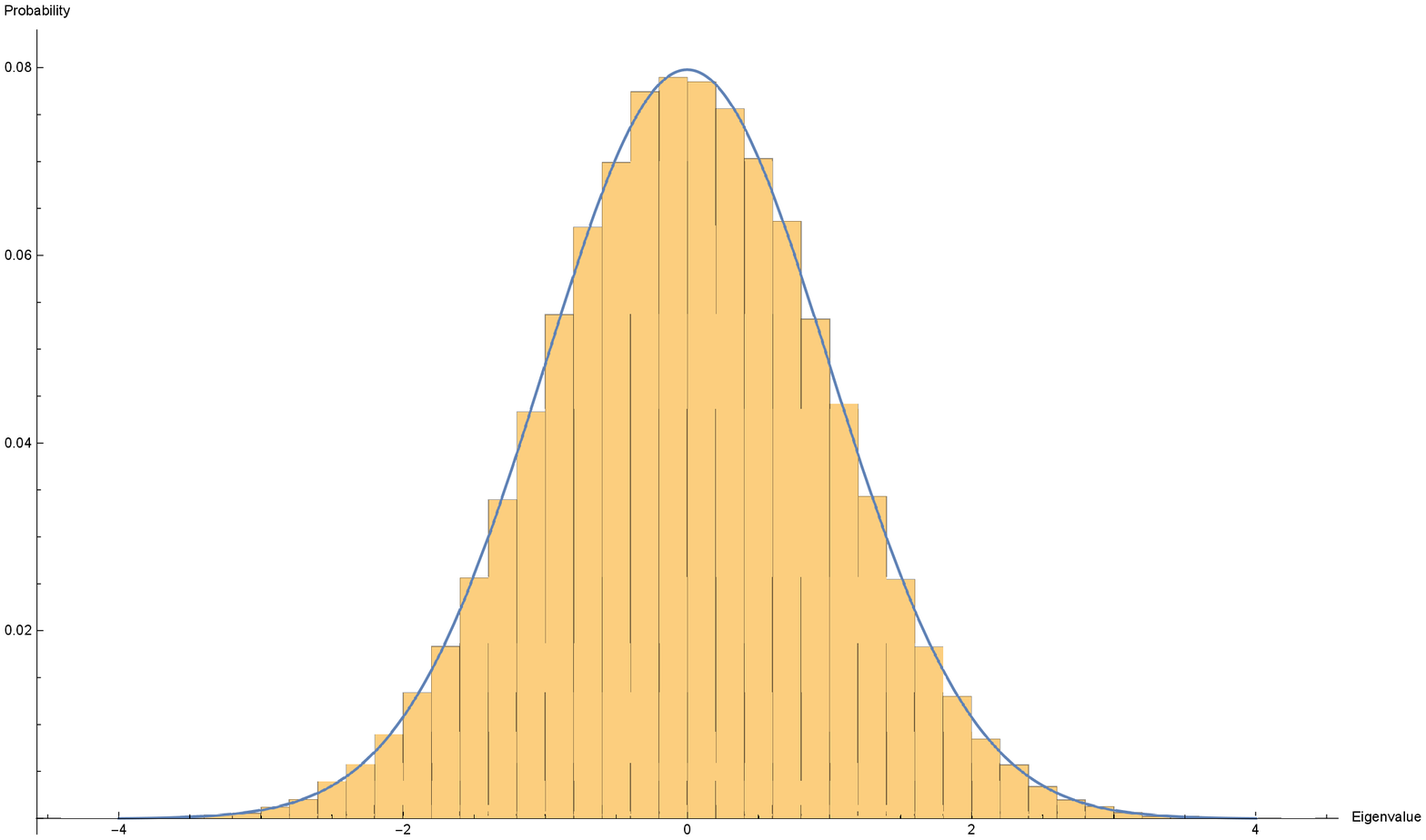}\ \ \  \includegraphics[scale=0.225]{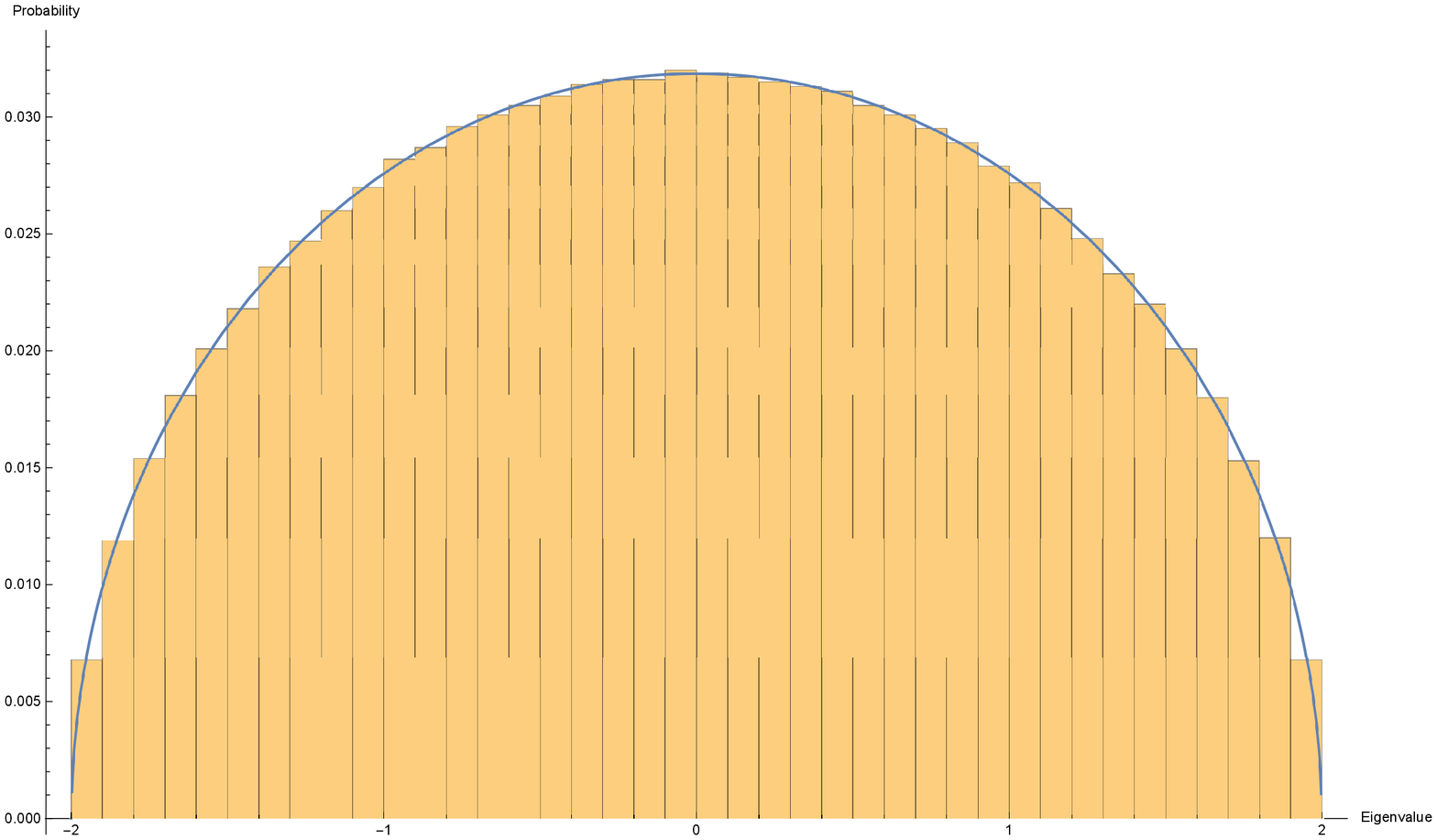}
     \caption{Eigenvalue distribution of $10,000 \times 10,000$ matrices: left is symmetric palindromic Toeplitz (plotted against a Gaussian), right is real symmetric (plotted against a semi-circular denisty).}
    \label{fig:my_label}
\end{figure}

The construction is of interest as a way to create ensembles and see how the properties of the constituent components are reflected in the new family; though inspired by a question from number theory, the resulting distributions do not correspond to those observed from convolving families of $L$-functions. Our analysis shows that the new
construction of \eqref{eq:1-disco-construction}
exhibits hybrid behaviors that bear resemblance to the
limiting distributions of its component matrices, converging to a new
universal distribution distinct from both the Gaussian and semicircular,
while retaining notable similarities to both. We then extend this construction
by considering
random block matrices constructed by successively
concatenating $\mathcal{D}_1 \left( A, B_1 \right)$ with additional matrices $\{B_k\}$.
Our work shows that
their scaled eigenvalues converge as $N \to \infty$.
An entire spectrum of fascinating hybrid behavior
exists
for the limiting eigenvalue distributions, uncovering
a galaxy of new, universal distributions.


\subsection{New Results}\label{subsec:Notation}

Entries of all matrices are defined on a  common probability space
$(\Omega,\mathcal F, \text{\rm P})$. Suppose $A$ and $\{B_k\}_{k \geq1}$ are  independent real symmetric random matrices, with possibly additional structure imposed. The matrix $A$ is of order $N\times N$ and $B_k$ is of order $2^kN\times 2^kN$. We shall assume that all the random entries have mean $0$ and variance $1$.

\begin{defi}For $d \in \mathbb{Z}^+$, the $d$-Disco of $A$ and $\mathbf{B} = \{B_k\}$,
denoted $\mathcal{D}_d (A, \mathbf{B} )$, is the $2^dN\times 2^dN$ real symmetric random matrix given by
\begin{equation}
    \mathcal{D}_d (A, \mathbf{B} )\ =\ \left[ \begin{array}{ccc}
        \begin{array}{cc}
         \begin{array}{cc}
            A & B_1 \\
            B_1 & A
        \end{array} & \text{\LARGE $B_2$}\\
        \text{\LARGE $B_2$} & \begin{array}{cc}
            A & B_1 \\
            B_1 & A
        \end{array} \end{array} & \cdots & \text{\fontsize{45}{0} $B_d$} \\
        \vdots & \ddots & \vdots \\
        \text{\fontsize{45}{0} $B_d$} & \cdots &  \begin{array}{cc}
         \begin{array}{cc}
            A & B_1 \\
            B_1 & A
        \end{array} & \text{\LARGE $B_2$}\\
        \text{\LARGE $B_2$} & \begin{array}{cc}
            A & B_1 \\
            B_1 & A
        \end{array} \end{array}
    \end{array} \right].
\end{equation}
\end{defi}
$\mathcal{D}_d$ depends on $A$ and $\{B_k\}$ but we suppress this dependence for simplicity.
We can write $\D_d$ inductively as
\[
    \D_d \ = \  \begin{pmatrix}
            \D_{d -1} & B_{d} \\
            B_{d} & \D_{d - 1}
            \end{pmatrix}.
\]

Observe that \eqref{eq:1-disco-construction} is a specific instance of the preceding construction.

\begin{defi}[Empirical spectral
distribution and measure]\label{defi:nesdab} Suppose $M_n$ is a real symmetric $n\times n$ random matrix with eigenvalues
 $\lambda_1, \ldots, \lambda_{n}$. The \emph{empirical spectral distribution} (ESD) is defined by
\begin{equation}
F_{M_n }(x) \ =\ \frac{\#\{i \le n:
\lambda_i\le x\}}{n}.
\end{equation}
The corresponding measure $\mu_{M_{n}}$ is called the empirical spectral measure.
The {\it expected empirical spectral distribution function} (EESD) of $M_n$ is defined as
$$\E(F_{M_n}(x))\ = \ \frac{1}{n}\sum_{i=1}^{n}\p [\lambda_i \leq x], \ x\in \mathbb{R}.$$
It is a non-random distribution function and we shall write $\E(F_{M_n})$ in short. The corresponding probability law  is known as the {\it expected spectral measure} of $M_n$.
\end{defi}

\begin{defi} \label{weak_conv} For any probability distribution function $F$ on $\mathbb{R}$, let
$$C_F \ = \  \{t: t \ \text{is a continuity point of}\  F\}.$$
A sequence of probability distribution functions $\{F_n\}$ is said to converge \textit{weakly} to a  probability distribution function $F$
if
\begin{equation*} F_n(t)\to F(t)\ \text{for all} \ t \in C_F.
\end{equation*}
\end{defi}
Our interest is in the convergence of the ESD and EESD of random real symmetric matrices as their dimension tends to $\infty$.

\begin{defi} \label{def:probasesd}
Let $F$ denote  a \textit{non-random} distribution function, defined on $\mathbb{R}$.

\begin{itemize}
    \item[(a)] the ESD of $M_n$ converges to $F$ \textit{almost surely} if for almost every $\omega \in \Omega$ and for all  $t\in C_F$,
$$F_{M_n}(t)\rightarrow F(t) \ \ \mbox{as}\ n\rightarrow \infty.$$

\item[(b)] The EESD of  $\{M_n\}$ converges to $F$ if $\E(F_{M_n})$ converges \textit{weakly} to $F$.

\end{itemize}

\end{defi}
It is easy to see that (a) $\Rightarrow$ (b).
We refer to the limit as the \textit{limiting spectral distribution} (LSD) of $\{M_n\}$ and the corresponding probability law as the \textit{limiting spectral measure}.

The following result is well-known\footnote{For a quick proof, see Lemma 1.2.1 of  \cite{B}.}.

\begin{lem}\label{lem:ch1.l3}
Suppose that $\{Y_n:n\in\bbn\}$ is a sequence of real-valued random variables with distribution functions $\{G_n\}$ such that for all $k\in\bbn$,
\[
\lim_{n\to\infty}\E(Y_n^k)\ = \ m_k\,\; \text{(finite)},
\]
and that there is a unique distribution function $G$ whose $k$-th moment is $m_k$ for every $k$. Then
$G_n\to G$.
\end{lem}

Suppose $M_n$ is real symmetric. Then the $h$-th moment of its ESD and EESD are given by
\begin{align}
m_h(F_{M_n})
&= \frac{1}{n} \sum_{i=1}^n
\lambda_i^h \ = \  \frac{1}{n} \mathrm{Tr}(M_n^h) \ = \  \mathrm{Tr} (M_n^h)\ = \ m_h (M_n) \label{eq:tracemoment}\\
m_h(\E(F_{M_{n}}))&= \ \E\left[\frac{1}{n} \mathrm{Tr}(M_n^h)\right]\ = \ \E \left[m_h(M_n)\right]\ \ .
\end{align}
where $\mathrm{Tr}$ denotes the trace.  (\ref{eq:tracemoment}) is known as the \textit{Trace-Moment} formula.
Now consider the following conditions.\vskip5pt

\begin{itemize}
    \item[(C1)] For every $h \geq 1$,  $\E[m_h (M_n)] \rightarrow m_h$.

    \item[(C2)] The moment sequence  $\{m_h\}$ corresponds to a unique probability distribution $F$.

    \item[(C3)] For every $h \geq 1$,  $$\sum_{n=1}^\infty[m_h (M_n)-\E[m_h (M_n)]^4 < \infty.$$
\end{itemize}

The next lemma follows easily from Lemma \ref{lem:ch1.l3}. We omit its proof.

\begin{lem}\label{lem:momenttrace}
Suppose $M_n$ is an $n \times n$ real symmetric matrix satisfying (C1), (C2). Then as $n \to \infty$ the EESD of $M_n$ converges weakly to $F$ determined by $\{m_h\}$. If $M_n$ also satisfies (C3), then the ESD converges weakly almost surely to $F$.
\end{lem}

To prove the convergence of the EESD and ESD of $\D_d$, it suffices to verify conditions (C1), (C2), and (C3). We record the assumption on the entries of our random matrices that will be called frequently in this article.
\vskip5pt

\noindent \textbf{Assumption 1}. Suppose the collection $\mathcal{C}$ of random variables from which the matrices are formed, are independent, have mean $0$, and variance $1$. Further,
$\sup_{X\in \mathcal{C}}\E|X|^k < \infty$ for all $k \geq 1$.
\vskip5pt

The following two results are known.

\begin{itemize}
    \item[i.] The ESD of $A/\sqrt{N}$ converges to the standard Gaussian law, almost surely.
		\item[ii.] The ESD of $B_k/\sqrt{2^kN}$ converges to the semi-circular law, almost surely.
\end{itemize}
For proofs  of the above two results under several alternate assumptions, including under Assumption A, see Theorems 2.1.3 and 2.4.2 of \cite{B} respectively.

Now  consider $\D_1$.   It can be shown that
\begin{equation}\label{eq1}
     \Tr (\D_1^k) \ = \  \dfrac{1}{2}\left[ \Tr (A + B_1)^k + \Tr (A - B_1)^k\right].
\end{equation}

We now state our main result; after looking at moments to see some general properties of the new construction and reviewing needed results from free probability in \S\ref{sec:preliminaries}, we prove a special case first in \S\ref{subsec:specialresult} and then the general case in \S\ref{subsec:mainresult}.

\begin{thm}\label{thm:D_d}
Suppose $A$ is a palindromic Toeplitz matrix and $\{B_k\}$ are Wigner matrices. Suppose these matrices are independent and the entries of each matrix are independent with mean zero and variance 1 and satisfy Assumption 1. Then
the EESD of $\D_d$ converges as $N\to \infty$. The LSD is the law $\mu_d$ of the self-adjoint variable
\[
    a_d\ :=\ \bigg(\frac{1}{\sqrt{2}}\bigg)^dg + \sum_{i = 1}^d \bigg(\frac{1}{\sqrt{2}}\bigg)^{i} s_i,
\]
where $g$ is a standard Gaussian variable, $s_1, \ldots s_d$ are standard free-Gaussians (standard semi-circular variables) and they are jointly free. Moreover, as $d\to \infty$, the probability law $\mu_d$ converges to the standard semi-circular law.
\end{thm}

\section{Preliminaries}\label{sec:preliminaries}

\subsection{Method of Moments}
Comparing the moments of the disco matrix to those of the Gaussian and semi-circular distributions offers insight into how the disco structure creates a hybrid of disparate limiting distributions. To proceed by the Method of Moments, we need to compute the expected traces of non-commutative, bivariate matrix polynomials resulting from the expansion of Equation \eqref{eq1}. Similarly, for general $\D_d$, we need to handle expected traces of non-commutative, multivariate matrix polynomials. The problem of computing expected traces may be reformulated as a combinatorics problem in which one must pair $2k$ points on the circumference of a circle with chords possibly subject to additional constraints. In the Gaussian case no constraints are placed, and all possible pairings of points on a circle contribute equally to the moment in the limit $N\to \infty$. In contrast, the semi-circular case of the real symmetric matrices has equal contribution from all pairings that have no crossings, while pairings with a crossing contribute zero in the limit $N \to \infty$.

Computing the moments of $\mathcal{D}_d$ represents a hybrid of these conditions: elements corresponding to real symmetric submatrices are paired with non-crossing chords $c_{RS}$, while the elements corresponding to PST submatrices are paired with chords $c_{PST}$ allowed to cross each other but not those in $c_{RS}$.

\begin{figure}[ht]
    \centering
    \includegraphics[width=0.2\textwidth]{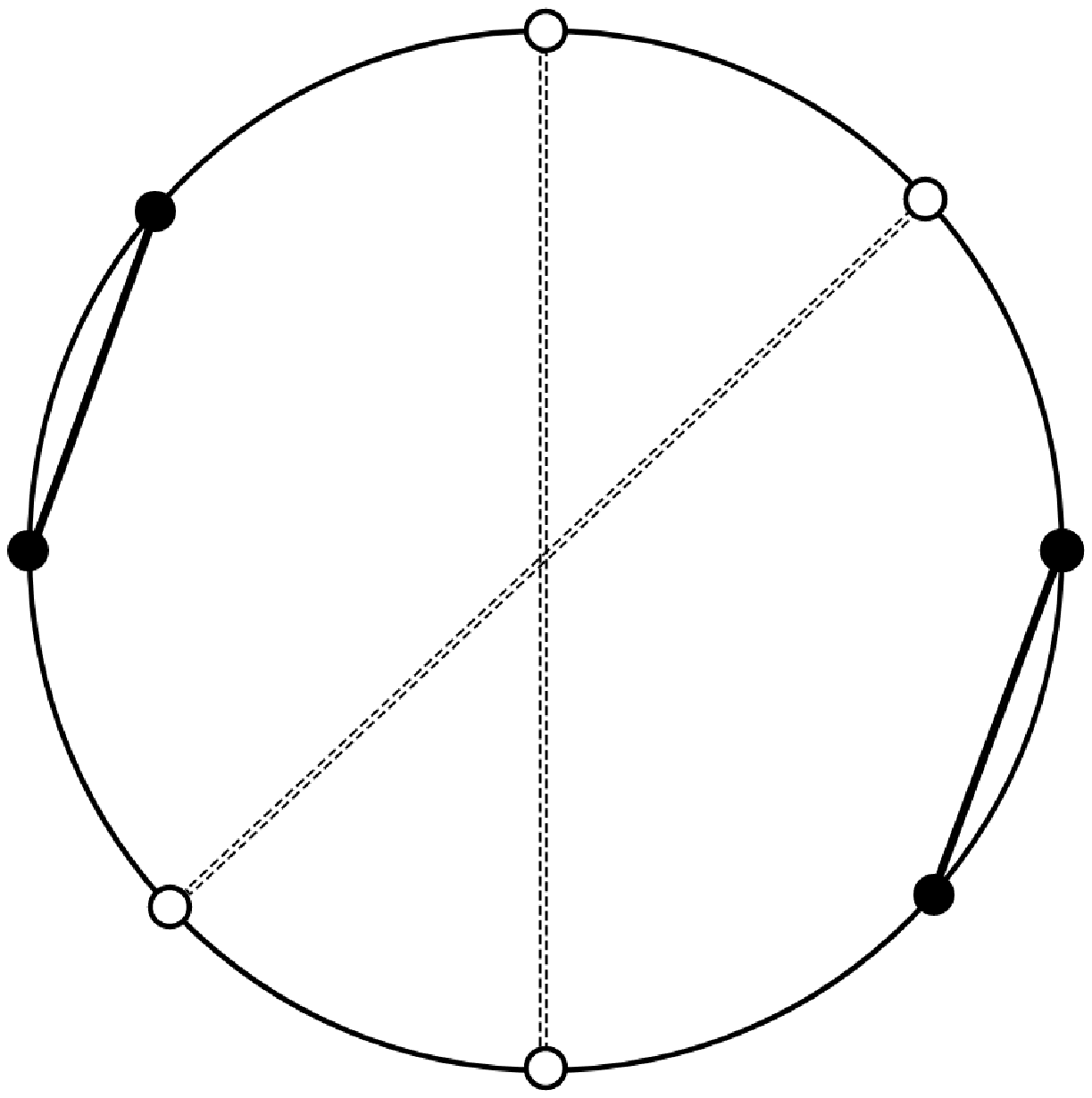}
    \includegraphics[width=0.2\textwidth]{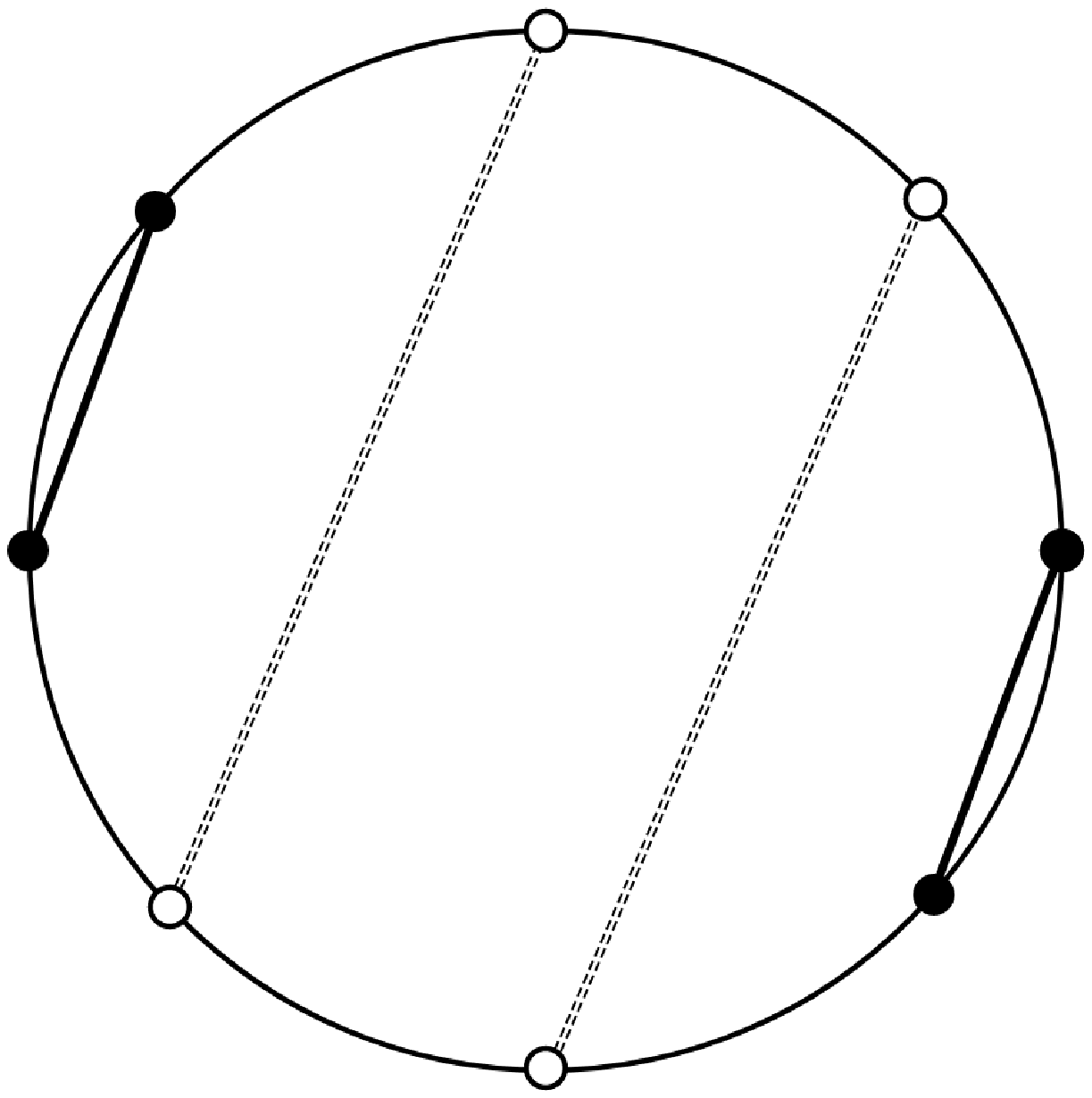}
    \includegraphics[width=0.2\textwidth]{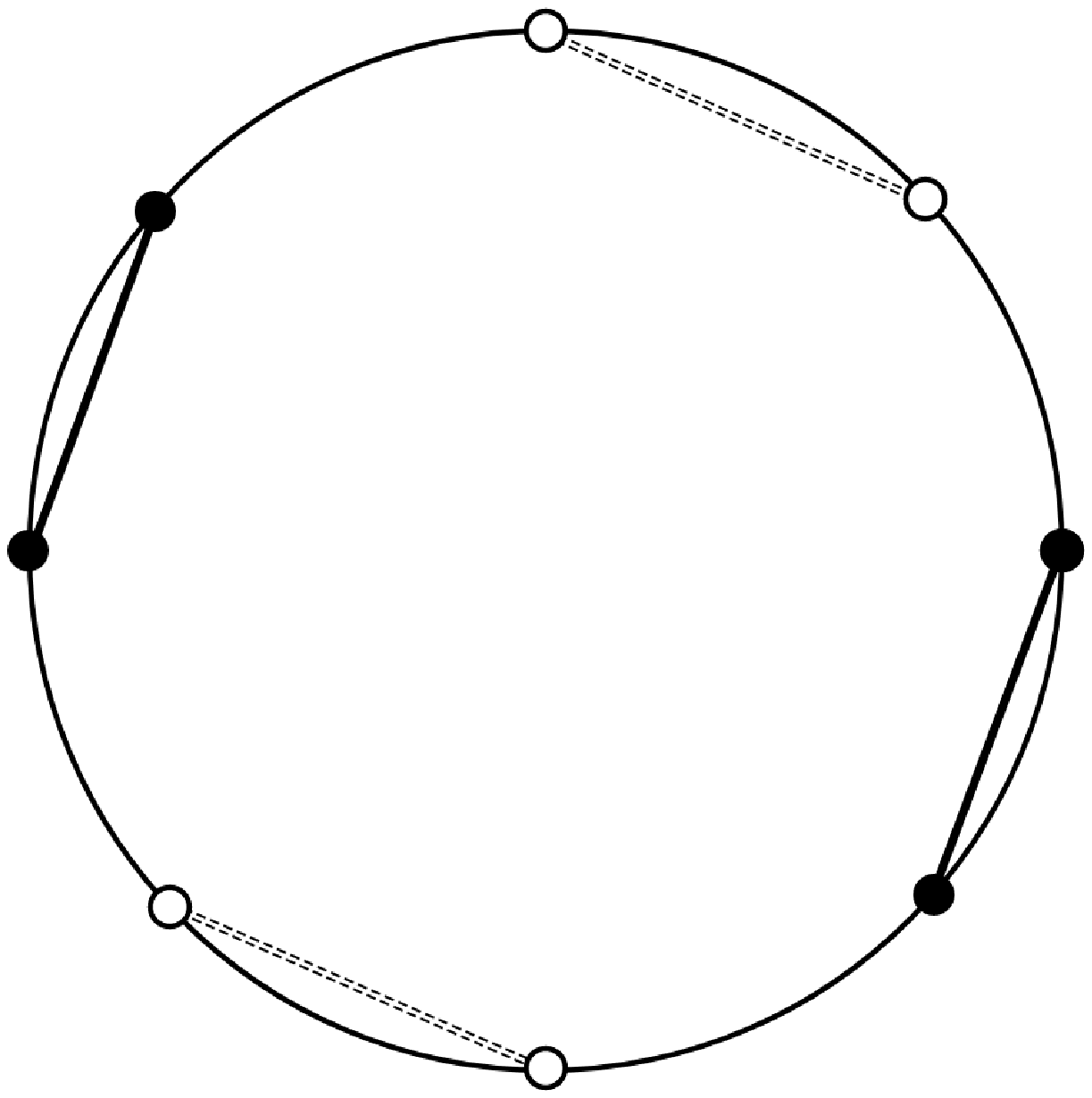}
    \includegraphics[width=0.2\textwidth]{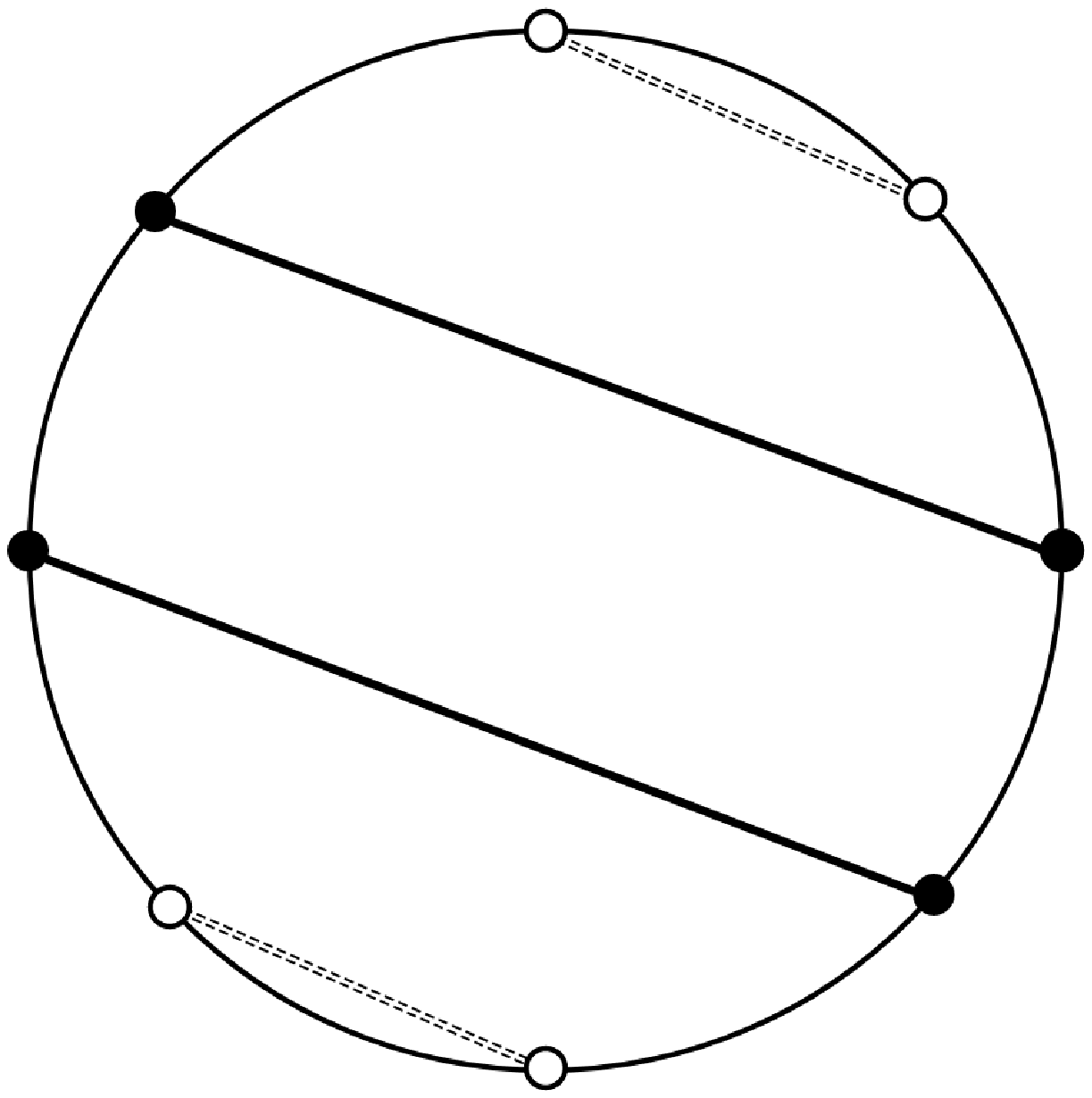}
    \caption{Visualization of contributing pairings of
    $\mathbb{E} \left[ A^2 B^2 A^2 B^2 \right]$.}
    \label{fig:contributing_pairings_A2B2A2B2}
\end{figure}

\begin{figure}[ht]
    \centering
    \includegraphics[width=0.2\textwidth]{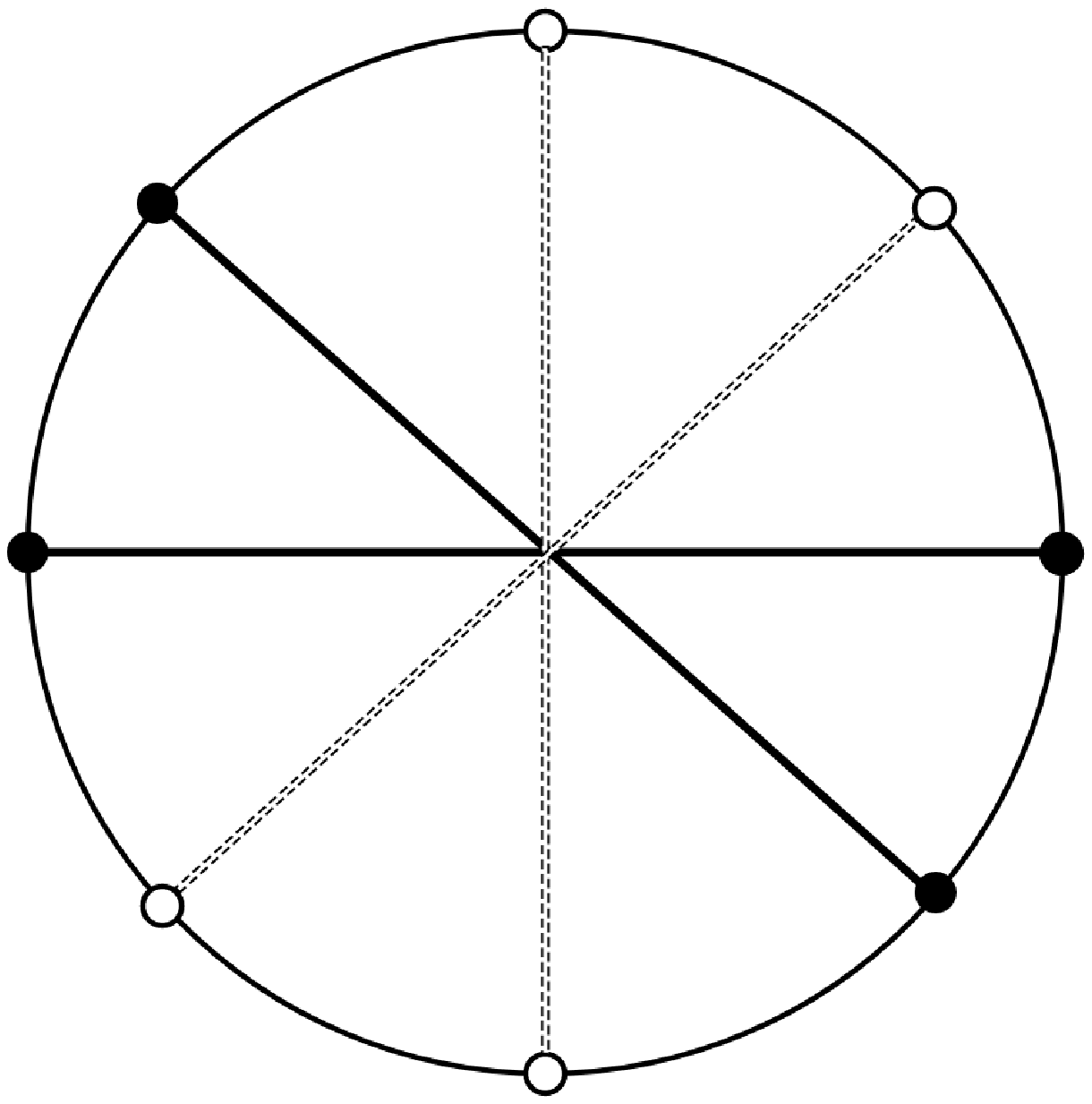} \includegraphics[width=0.2\textwidth]{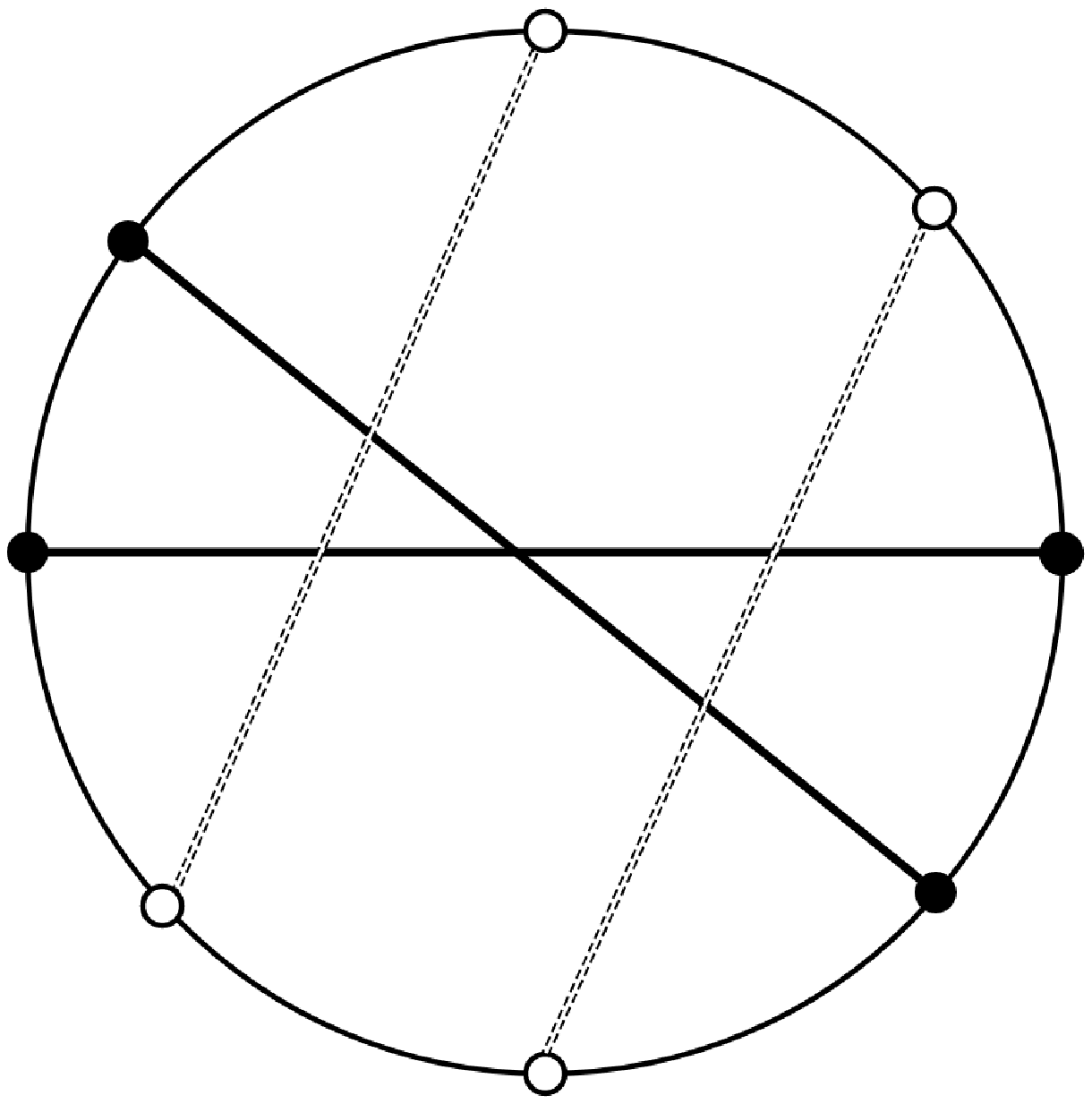}
    \includegraphics[width=0.2\textwidth]{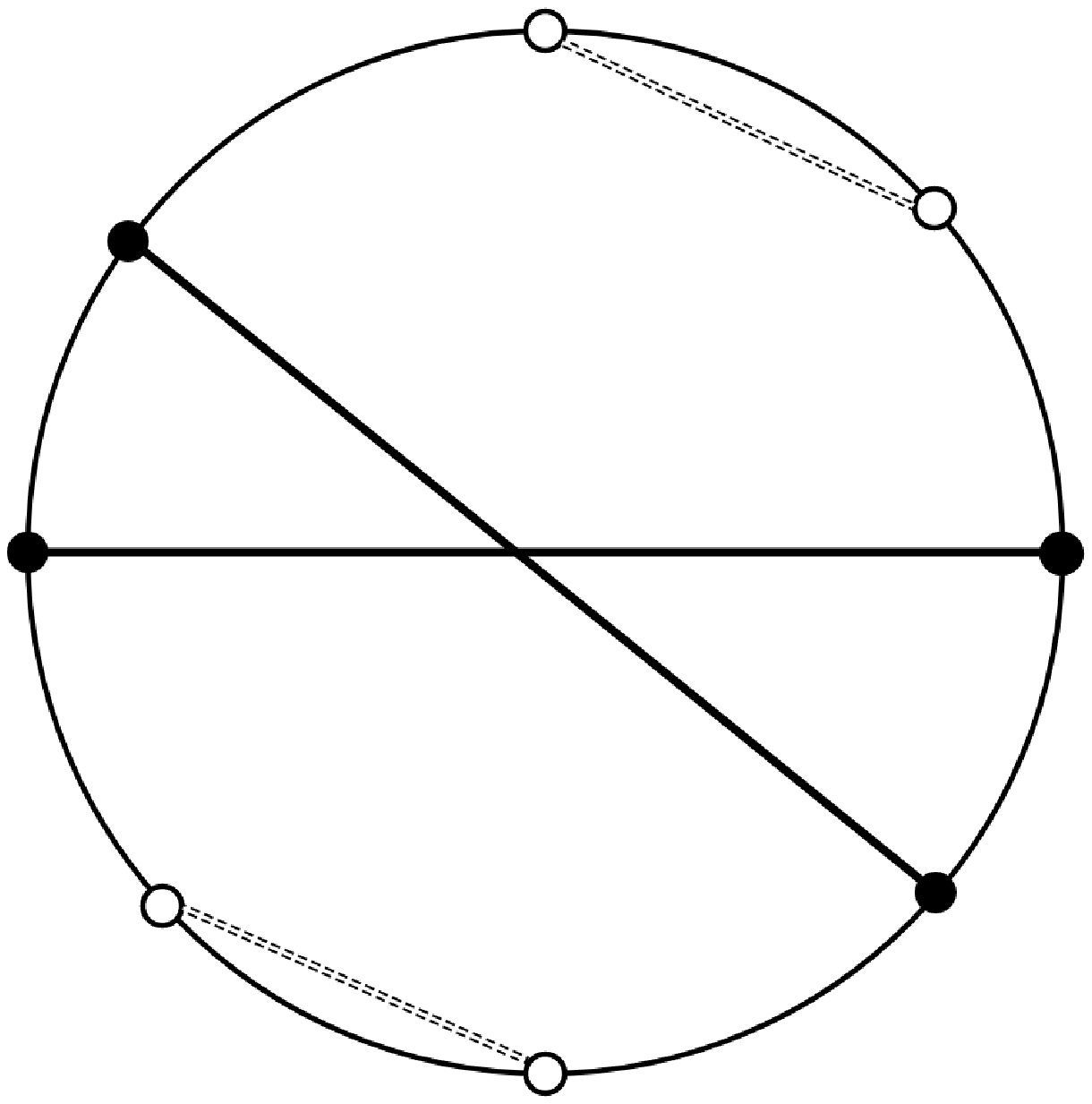}
    \includegraphics[width=0.2\textwidth]{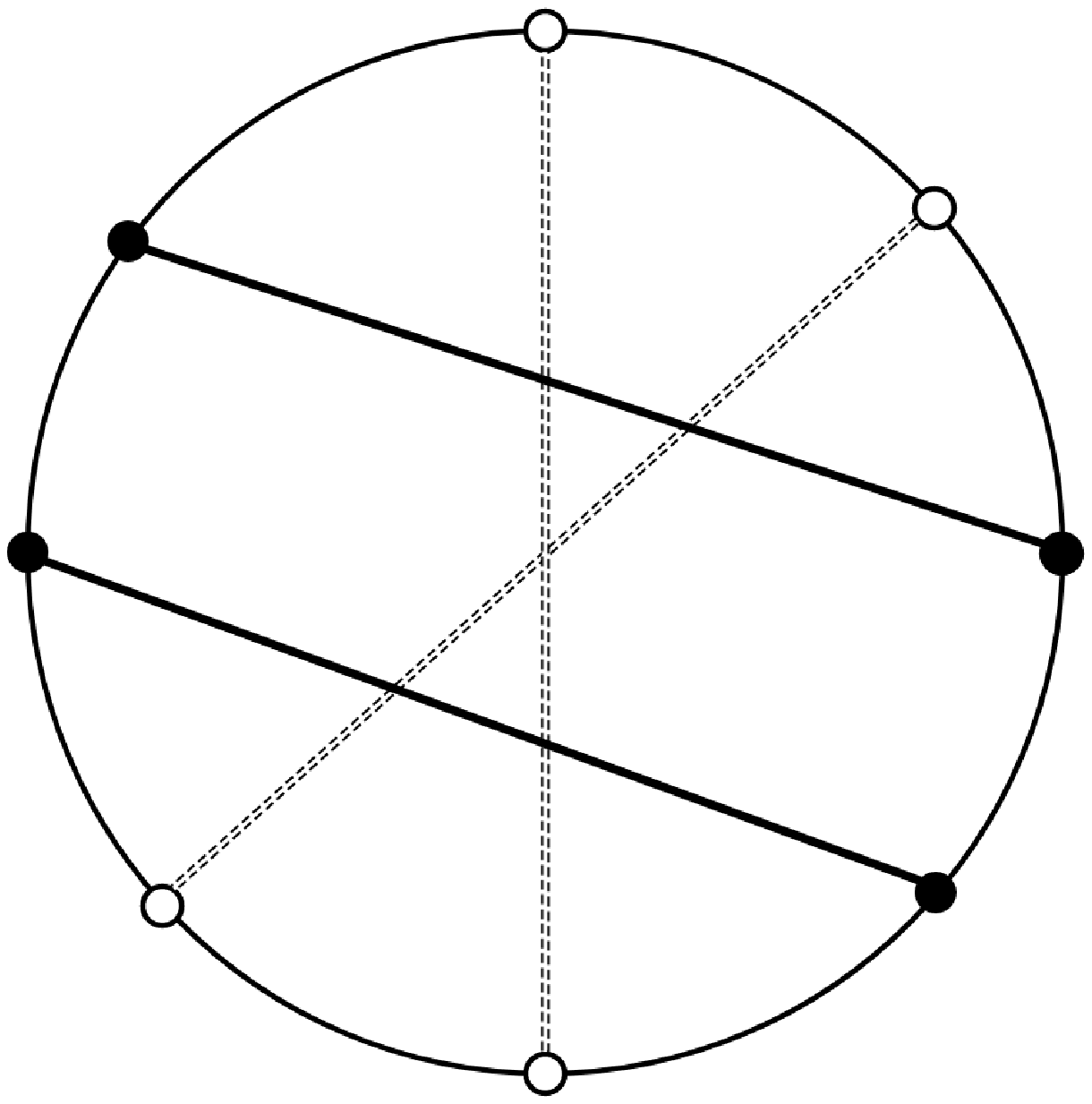}
    \caption{Visualization of non-contributing pairings of
    $\mathbb{E} \left[ A^2 B^2 A^2 B^2 \right]$.}
    \label{fig:non_contributing_pairings_A2B2A2B2}
\end{figure}

It readily follows that the odd moments of $\mathcal{D}_d$ are zero, while  brute-force computation of low even moments of $\mathcal{D}_1$ is consistent with hybrid distribution behavior.

\begin{center}
\begin{tabular}{|c|c|c|c|}\hline
Moment & Semi-circular & \hspace{0.5cm} $\mathcal{D}_1$ \hspace{0.5cm}
& Gaussian \\ \hline
2 & 1  &\ 1.00  & 1 \\ \hline
4 & 2  &\ 2.25  & 3 \\ \hline
6 & 5  &\ 7.00    & 15  \\ \hline
8 & 14 & 27.50 & 105 \\ \hline
\end{tabular}
\end{center}

\begin{figure}[ht]
    \centering
    \includegraphics[scale=0.27]{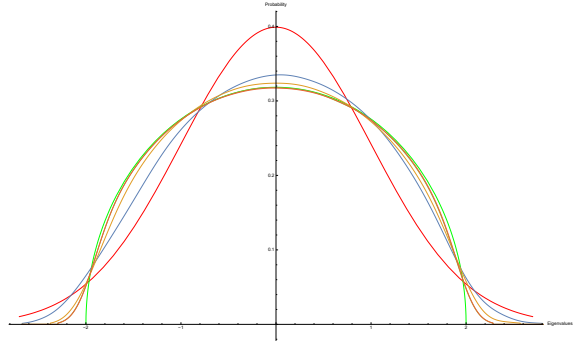}
    \caption{Eigenvalue distribution of $2^{14} \times 2^{14}$ matrices $\mathcal{D}_d$, $d = 1, 2, 4, 8$, alongside the Gaussian (red) and semi-circular (green).}
    \label{fig:disco-illustration}
\end{figure}

Elementary generalizations of such combinatorial analysis shows that the even moments of $\mathcal{D}_d$ are bounded between those of the semi-circular and Gaussian. For $d$ finite, the scaled eigenvalues converge as $N \to \infty$ to intermediate distributions with hybrid Gaussian and semi-circular behavior. Taking $d \to \infty$ causes the scaled eigenvalues to converge exponentially to the semi-circular distribution; see Figure \ref{fig:disco-illustration}. For a detailed proof of these results via the Method of Moments, see \cite{BBDLMMWX}.

On the other hand,  the question of convergence of $\D_d$ can be framed in terms of convergence of elements in $*$-probability spaces. The concept of (asymptotic) free independence plays a crucial role. Arguments based on this theory of non-commutative probability spaces, besides shortening the proof significantly, also provides better insight into the nature of the limiting spectral distribution of $\mathcal{D}_d$. We now provide a brief description of the essential background and results from this theory that we shall need for the proofs in \S\ref{sec:proofs}.

\subsection{Elements of free probability}


\begin{defi}[$*$-algebra]
\label{def:$*$-algebra} A collection $\mathcal{A}$ is called a unital algebra if it is an additive vector space over $\mathbb{C}$, endowed with multiplication satisfying the following for all $x, y, z\in \mathcal{A}$ and $\alpha\in \mathbb{C}$:  \vskip3pt

\begin{itemize}
    \item[(i)] $x(yz)=(xy)z$,

    \item[(ii)] $(x+y)z=xz+yz$,

    \item[(iii)] $x(y+z)=xy+xz$,

    \item[(iv)] $\alpha(xy)=(\alpha x)y=x(\alpha y)$,

    \item[(v)] there exists a multiplicative identity element $\one_{\mathcal{A}}$ in $\mathcal{A}$.
\end{itemize}

$\mathcal{A}$ is called a  \textit{$*$-algebra} if there exists a mapping $\mathcal{A} \to \mathcal{A}: x\mapsto x^*$ such that,
for all $x, y\in \mathcal{A}$ and $\alpha\in \mathbb{C}$:

\begin{itemize}

\item[(vi)] $(x+y)^*=x^*+y^*$,

\noindent
\item[(vii)] $(\alpha x)^*=\bar\alpha x^*$,

\noindent
\item[(viii)] $(xy)^*=y^*x^*$, and

\noindent
\item[(ix)] $(x^*)^*=x$.

\end{itemize}
\end{defi}

\begin{defi}[$*$-probability space] Suppose $\mathcal{A}$  is a unital $*$-algebra over $\mathbb{C}$ and $\varphi$ is a linear functional on $\mathcal{A}$  such that $\varphi$ is positive ($\varphi(aa^*)\geq 0$,
$\forall \; a\in \A$) and  $\varphi(\one_{\mathcal{A}})=1$. Then the pair $(\mathcal{A}, \varphi)$ is called a $*$\textit{-probability space}.
\end{defi}
\begin{exa}
Fix a positive integer $n$.
Let $\mathcal{M}_{n}(\mathbb{C})$ be the algebra of $n\times n$ matrices with complex entries under ordinary addition and multiplication. Consider the expected normalized trace:
\begin{equation}\E {\rm tr}(a)  \ = \  \dfrac{1}{n}\E \left[ {\rm Tr}(a) \right]\ = \ \dfrac{1}{d}\sum_{i=1}^{d}\E(\alpha_{ii})\, \ \ \forall\ \  a \ = \  ((\alpha_{ij}))_{i,j=1}^{d} \in \mathcal{M}_{n}(\mathbb{C}).
\nonumber
\end{equation}
Then $(\mathcal{M}_{n}(\mathbb{C}), \E{\rm{tr}})$ is a $*$-probability space.
\end{exa}

\begin{defi}[Probability laws of self-adjoint elements]\label{def:problaw} Suppose $(\mathcal{A}, \varphi)$ is a $*$-probability space and $a\in \A$ is a self-adjoint element. Then $\varphi(a^k), k\geq 1$, are called the \textit{moments} of $a$. If there is a probability law with moments $\varphi(a^k), k\geq 1$, and it is unique, then it is called the probability law of $a$ and is denoted by $\mu_a$.
\end{defi}

\begin{defi}[Semi-circular, circular and Gaussian elements]\label{def:freegauss}
Let $s$, $g$, and $c$ be self-adjoint elements.
\begin{itemize}
\item[(a)] $s$ is said to be (standard) semi-circular or free-Gaussian if its moments are given by
\begin{eqnarray} \label{eqn: semimoment4.1}
\varphi(s^h) \ = \  \begin{cases} C_n\ \ \text{{\rm if}\ $h=2n$}, \\
0\ \ \text{{\rm if}\ $h$ {\rm is odd},}
\end{cases}
\end{eqnarray}
where $\{C_n\}$ are the Catalan numbers
$$C_n\ = \ \frac{1}{n+1}\binom{2n}{n},\ \ n\geq 1.$$
Then $\mu_s$ is called the (standard) semi-circular law.
\vskip3pt

\item[(b)] $g$ is said to be (standard) Gaussian if its moments agree with the Gaussian moments. Its probability law will be denoted by $\mu_g$.

\vskip3pt
\item\item[(i)]  $c$ is said to be circular if its moments are given by:
\begin{eqnarray}\label{eqn:circularmoments}
\varphi(c^{\epsilon_1}c^{\epsilon_2}\cdots c^{\epsilon_p}) \ = \  \begin{cases}
\sum\limits_{\pi=\{\{(r_i,s_i)\}}\displaystyle{\prod_{i=1}^k}
\mathbb{I}\{\epsilon_{r_{i}} \neq \epsilon_{s_{i}}\} \ \ \text{if}\ \  p=2k,\\
0\ \ \text{if}\ \ p=2k+1
\end{cases}
 \end{eqnarray}
for all choices of $\{\epsilon_i\}$ from the set $\{1, *\}$. The sum is over all non-crossing pair-partitions $\pi=\{(r_i, s_i), r_i < s_i, 1\leq i \leq k\}$ of $\{1, \ldots, 2k\}$.
\end{itemize}
\end{defi}

The concept of free independence for $*$-probability spaces is the analogue of classical independence for probability spaces; for more details, see \cite{NS}.


\begin{defi}[Free independence]\label{def:freemoments}
For a fixed index set $I$,
 let $(\A_i,  i\in I)$ be unital $*$-subalgebras of $(\A, \varphi)$. Then
$(\A_i, i\in I)$ are called \textit{freely independent}, if $\varphi(a_1 \cdots a_k) \ = \  0$
for every $k\geq 1$ and for every $1\leq j\leq k$ whenever,
(a) $a_j\in\A_{i(j)}$ ($i(j)\in I$),  (b) $\varphi(a_j)=0$ for every $1\leq j\leq k$ and, (c)  neighboring elements are from different sub-algebras, that is, $i(1)\neq i(2)\neq \cdots i(k-1)\neq i(k)$.
Elements are called freely independent, or simply, free, if the $*$-sub-algebras generated by them are free.
\end{defi}

It can be shown that if $s_1$ and $s_2$ are two free standard semi-circular elements, then $(s_1+\sqrt{-1} s_2)/\sqrt{2}$ has the same moments as a circular variable.

\begin{defi}[Free additive convolution]\label{def:freeconv}  Suppose $a$ and $b$ are two self-adjoint elements from a $*$-probability space  with probability laws $\mu_a$ and $\mu_b$ respectively. Then if $\mu_{a+b}$ exists, it is called the free additive convolution of $\mu_a$ and $\mu_b$.
\end{defi}

We now explain the notion of convergence for elements in $*$-probability spaces. We shall use the  notation $\Pi (\cdot)$ to denote polynomials formed from elements of a $*$-algebra. By default, adjoints are included in the arguments of the polynomials.

\begin{defi}[Convergence of elements in $*$-probability spaces] \label{def: convergence}
Let $(\mathcal{A}_{n}, \varphi_n)$, $n \geq 1$ be a sequence of $*$-probability spaces and let $(\mathcal{A}, \varphi)$ be another $*$-probability space.
\vskip3pt

\begin{itemize}
\item[(a)] (Marginal convergence) We say that $a^{(n)} \in \mathcal{A}_{n}$
\textit{converges in $*$-distribution} to $a \in \mathcal{A}$ if
\begin{equation} \label{eqn: arr4.11}
\lim\varphi_{n}(\Pi(a^{(n)}))\ = \  \varphi(\Pi(a_i)) \ \ {\rm for \ all \ polynomials}\ \ \Pi. 
\end{equation}
We denote this convergence by $a^{(n)}\stackrel{*}{\rightarrow} a$.
\vskip3pt

\item[(b)] (Joint convergence)\index{joint convergence} Suppose $I$ is an index set. The elements
$\{a_i^{(n)}: i \in I\}$ from $\mathcal{A}_{n}$ are said to \textit{converge (jointly)} to $\{a_i: i \in I\}$ from $\mathcal{A}$ if,
\begin{equation}\label{eqn:ggrt4.1}
\Pi(\{a_i^{(n)}: i \in I\}) \stackrel{*}{\rightarrow} \Pi(\{a_i: i \in I\}) \ \ \text{for all polynomials} \ \
\Pi.\end{equation}
We write this as $\{a_i^{(n)}: i \in I\}\ \stackrel{*}{\rightarrow} \{a_i: i \in I \}$.
\vskip3pt

\item[(c)] (Asymptotic freeness) If  the limit variables $\{a_i, i\in I\}$ are free, then we say that
$\{a_i^{(n)}, i \in I\}$  are free in the limit, or, are \textit{asymptotically free}.

\end{itemize}
\end{defi}

\begin{exa} Consider the $*$-probability spaces $(\mathcal{M}_n(\mathbb{C}), \E{\rm tr}), n \geq 1$. Suppose $M_n\in \mathcal{M}_n(\mathbb{C}), \ n \geq 1$ are real symmetric. Then $M_n$ converges in $*$-distribution if and only if
\begin{equation}\varphi_n(M_n^k)\ = \ \frac{1}{n}\E \left[ {\rm Tr}(M_n^k) \right] \ \ {\rm converges \ for \ all}\ \ k \in \bbn.\nonumber
\end{equation}
The limit algebra is generated by an indeterminate $m$, say, and $\varphi$ is defined on this algebra by
\begin{equation}\varphi (m^k)\ = \ \lim_{n\to \infty} \frac{1}{n} \E \left[ {\rm Tr}(M_n^k) \right]\ = \  m_k\ \ {\rm (say)}.\nonumber\end{equation}
Thus the convergence in $*$-distribution of $\{M_n\}$ is the same as the convergence of  $m_k(\E F_{M_n})$ for $k\geq 1$.

Note that there is always a probability law, say $\mu$, with $m_k, k\geq 1$ as its moments. If  $\mu$ is unique, then EESD of $M_n$ converges weakly to $\mu$ by Lemma \ref{lem:ch1.l3}. This connects the convergence of EESD to convergence in $*$-distribution.

If we have more than one sequence of matrices, then, likewise, their joint convergence in $*$-distribution (with respect to $\E {\rm tr}$) is the same as the convergence of traces of all monomials in these matrices and their adjoints. If the matrices under consideration are real symmetric, then adjoints are redundant.
\end{exa}

We summarize two well-known\footnote{For the proofs of parts (a) and (b), see \cite{bose2011convergence}, and \cite{adhikari2019brown} respectively.} facts on $*$-distribution convergence.

\vskip5pt
\begin{cla}\label{cla:star-conv}
Let $W_i, 1\leq i \leq d$ be independent $n \times n$ Wigner matrices, $C_i, 1\leq i \leq d$ be independent $n \times n$ random matrices, and $A$ be a random PST matrix, where all entries satisfy Assumption 1. Then,

\begin{itemize}
\item[(a)]
    $\{A/\sqrt{n}, W_1/\sqrt{n}, \ldots , W_d/\sqrt{n}\}$ converge jointly in $*$-distribution as elements of $(\mathcal{M}_n(\mathbb{C}), \E {\rm tr})$ to $\{g, s_1, s_2, \ldots , s_d\}$ where $g, s_1, \ldots,  s_d$ are free, $g$ is Gaussian and $s_i$ are semi-circular variables in some $*$-probability space $(\A, \varphi)$.

\item[(b)]
    $\{C_1/\sqrt{n}, \ldots , C_d/\sqrt{n}, W_1/\sqrt{n}, \ldots , W_d/\sqrt{n}\}$ converge jointly in $*$-distribution as elements of $(\mathcal{M}_n(\mathbb{C}), \E {\rm tr})$ to $(c_1, \ldots c_d, s_1, \ldots, s_d)$ which are all free, and  $\{c_i\}$ are circular elements while $\{s_i\}$ are semi-circular variables.

\end{itemize}
\end{cla}

\section{Proofs}\label{sec:proofs}

\subsection{Convergence of \texorpdfstring{$\D_1$}{Lg}}\label{subsec:specialresult}
Observe that  we may diagonalize
$\mathcal{D}_1$ in the following manner:
\begin{align}
     \  \begin{bmatrix} I/\sqrt{2} & I/\sqrt{2} \\
                                 I/\sqrt{2} & -I/\sqrt{2}
							     \end{bmatrix}
									              \mathcal{D}_1
		             \begin{bmatrix} I/\sqrt{2} & I/\sqrt{2} \\
                                 I/\sqrt{2} & -I/\sqrt{2}
							     \end{bmatrix}
				&\ = \      		\begin{bmatrix} A+B_1 &   0 \\
                                  0   &  A-B_1 \end{bmatrix}
		 \end{align}
By Claim \ref{cla:star-conv}, $(A/\sqrt{N}, B_1/\sqrt{N})\stackrel{*}{\to} (g, s_1)$;
  it follows that
\begin{eqnarray*}\lim \E {\rm tr}\left(\frac{\D_1}{\sqrt{2N}}\right)^k  &\ = \ &\lim \dfrac{1}{2N}\E \left[ {\rm  Tr}\left(\frac{\D_1}{\sqrt{2N}}\right)^k \right]\\ &\ = \ & \lim\dfrac{1}{2}\left[\E \left[{\rm Tr}\left(\frac{A + B_1}{\sqrt{2N}}\right)^k \right] +  \E \left[{\rm Tr}\left(\frac{A - B_1}{\sqrt{2N}}\right)^k \right]\right]\\
                                  &\ = \ & \dfrac{1}{2}\left[\varphi\left(\frac{g + s_1}{\sqrt{2}}\right)^k +
																	\varphi\left(\frac{g - s_1}{\sqrt{2}}\right)^k\right]\\
&\ = \ & \varphi\left(\frac{g + s_1}{\sqrt{2}}\right)^k.
\end{eqnarray*}
The last equality follows since $(g,s_1)$ and $(g, -s_1)$ have the same moments.
In particular, the moments of $\D_1$ converge as $N\to \infty$. Thus $\{\varphi(\frac{g+s_1}{\sqrt{2}})^k)\}$  are the moments of a unique probability law, as $g$ and $s_1$ are free, have the Gaussian and the semi-circular laws respectively, and their moments are upper bounded\footnote{For a direct proof of this upper boundedness without the use of free probability, see \cite{BBDLMMWX}.} by Gaussian moments. Indeed, this law is the free additive convolution of the probability laws $\mu_{s_1/\sqrt{2}}$ and $\mu_{g/\sqrt{2}}$, and is written as
$\mu_{g/\sqrt{2}}\boxplus \mu_{s_1/\sqrt{2}}$.
Thus we have arrived at the following theorem.

\begin{thm}\label{thm:d1}
Suppose $A$ is a palindromic Toeplitz matrix and $B_1$ is a Wigner matrix independent of $A$. Suppose the entries of each matrix are independent with mean zero and variance 1 and satisfy Assumption 1.
Then the EESD of $\D_1$ converges weakly. The LSD is the law, say $\mu_1$,  of the self-adjoint variable $a_1=(g+s)/\sqrt{2}$ and is  the free additive convolution of $\mu_{g/\sqrt{2}}$ and
$\mu_{s/\sqrt{2}}$.
\end{thm}

\begin{rek}\label{rek:noniid}
    By verifying (C3), it can be easily shown\footnote{See \cite{bose2011convergence} for similar arguments.} that the ESD converges almost surely to the same law; we omit the details.
\end{rek}


\subsection{Convergence of \texorpdfstring{$\D_d$}{Lg}}\label{subsec:mainresult}
Extending the ideas of the previous section yields our main result, which we restate for the convenience of the reader. \\ \

\noindent \textbf{Theorem \ref{thm:D_d}:} \emph{Suppose $A$ is a palindromic Toeplitz matrix and $\{B_k\}$ are Wigner matrices. Suppose these matrices are independent and the entries of each matrix are independent with mean zero and variance 1 and satisfy Assumption 1.
Then the EESD of $\D_d$ converges as $N\to \infty$. The LSD is the law $\mu_d$ of the self-adjoint variable}
\[
    a_d\ :=\ \bigg(\frac{1}{\sqrt{2}}\bigg)^dg + \sum_{i = 1}^d \bigg(\frac{1}{\sqrt{2}}\bigg)^{i} s_i,
\]
\emph{where $g$ is a standard Gaussian variable, $s_1, \ldots s_d$ are standard free-Gaussians (standard semi-circular variables) and they are jointly free. Moreover, as $d\to \infty$, the probability law $\mu_d$ converges to the standard semi-circular law.} \\ \

Remark \ref{rek:noniid} also holds for $\mathcal{D}_d$, \emph{mutatis mutandis}. It may be noted that $\sum_{i = 1}^d \left(1/\sqrt{2}\right)^{i} s_i$ is a semi-circular variable with mean $0$ and variance $\sum_{i = 1}^d \left(\frac{1}{2}\right)^{i} = 1 + o_{d \rightarrow \infty}(1)$.

\begin{proof}Let us first consider $\D_2$. Note that
\[
    \frac{\D_2}{\sqrt{4N}} \ = \  \frac{1}{\sqrt{2}} \begin{pmatrix}
        \frac{1}{\sqrt{2N}}\D_1 & \frac{1}{\sqrt{2N}} B_2 \\
        \frac{1}{\sqrt{2N}} B_2 & \frac{1}{\sqrt{2N}}\D_1
    \end{pmatrix}.
\]
We have previously shown the following:
\begin{itemize}
\item[(i)] \ $\D_1$ is an element of $\mathcal{M}_{2N}(\mathbb{C})$ and $\D_1/\sqrt{2N} \stackrel{*} \to a_1$; and

\item[(ii)] $B_2$ is an element of $\mathcal{M}_{2N}(\mathbb{C})$  and $B_2/\sqrt{2N} \stackrel{*}\to s_2$, where $s_2$ is semi-circular.
\end{itemize}

We now claim that $\D_1/\sqrt{2N}$ and $B_2/\sqrt{2N}$ are asymptotically free.
To see this,  first partition $B_2$ as
\[
    B_2 \ = \  \begin{pmatrix}
        V & U\\
        U^\top & W
    \end{pmatrix},
\]
where $V$, $W$ are $N \times N$ Wigner matrices and $U$ is an $N \times N$ random matrix; the matrices $V$, $W$ and $U$ are independent. Note also that $U, V, W$ are independent of $A$ and $B_1$. Observe:

\begin{itemize}
\item[(iii)] From the results of \cite{basu2012joint},
$A$, $B_1$ $V$, and $W$ converge jointly and are asymptotically free.

\item[(iv)] From the results of \cite{adhikari2019brown} $U$, $V$, $W$ and $B_1$ converge jointly and are asymptotically free.
\end{itemize}

We omit a straightforward extension of the arguments in \cite{basu2012joint} and \cite{adhikari2019brown} proving that, as elements of $\mathcal{M}_{N}(\mathbb{C})$, $(A, U, V, W, B_1)\stackrel{*}{\to} (g, u, v, w, s_1)$ where these five variables are free, $g$ is standard Gaussian, $u$ is  circular and, $v, w, s_1$ are standard  semi-circular. Hence, using a decomposition like \eqref{eq1} on $\D_2$, and arguing as before, we may conclude that
\[
    \frac{\D_2}{\sqrt{4N}} \stackrel{*}\ \to\ \gamma_2\ =\ \frac{1}{\sqrt{2}} \begin{pmatrix}
    \frac{1}{\sqrt{2}}\begin{pmatrix}g & s \\ s & g \end{pmatrix} & \frac{1}{\sqrt{2}}\begin{pmatrix}v & u \\ u^* & w \end{pmatrix} \\
    \frac{1}{\sqrt{2}}\begin{pmatrix}v & u \\ u^* & w \end{pmatrix} & \frac{1}{\sqrt{2}}\begin{pmatrix}g & s \\ s & g \end{pmatrix}
\end{pmatrix}\ \stackrel{*}{=}\ \frac{g}{(\sqrt{2})^2} + \frac{s_1}{\sqrt{2}} + \frac{s_2}{(\sqrt{2})^2},
\]
where $g$ is standard Gaussian, $s_1$ and $s_2$ are standard semi-circular, and $g, s_1, s_2$ are free. The EESD of $\D_2/\sqrt{4N}$ converges to the distribution of $\mu_{g/(\sqrt{2})^2} \boxplus \mu_{s_1/\sqrt{2}} \boxplus \mu_{s_1/(\sqrt{2})^2}=\mu_2$. The case of general $d$ may be tackled similarly by an induction argument.




Finally, it is known that sum of free semi-circular variables is again semi-circular (the variances are summed). Thus, from the representation of $a_d$, we have a negligible term $g/(\sqrt{2})^d$ plus a free semi-circular variable whose variance is $1 + o_{d \rightarrow \infty}(1)$. Therefore, as $d \to \infty$, the LSD of $\mathcal{D}_d/\sqrt{2^d N}$, i.e. the law of $a_d$, approaches the semi-circular/free-Gaussian law.
\end{proof}

\section{Future Work}\label{sec:futurework}

So far we have investigated the density of the eigenvalues; we now
consider another problem, that of the spacings between adjacent
eigenvalues of $\mathcal{D}_d \left( A, \mathbf{B} \right)$. For
concreteness we restrict our purview to $A$ a symmetric palindromic
Toeplitz matrix and $\mathbf{B}$ a sequence of real symmetric
matrices. In this case $\mathcal{D}_d$ has only
\begin{equation}
    \frac{N}{2} + \sum_{i=1}^d \frac{2^{i-1} N (2^{i-1} N + 1)}{2}
\end{equation}
degrees of freedom, which is much smaller than $2^d N (2^d N + 1)/2$,
it is reasonable to believe the spacings between adjacent normalized
eigenvalues
$(\gl_{i+1}(\mathcal{D}_d) - \gl_i(\mathcal{D}_d)) / \sqrt{N}$
may differ from those of full real symmetric matrices.

\begin{figure}[ht]
    \centering
    \includegraphics[width=0.6\textwidth]{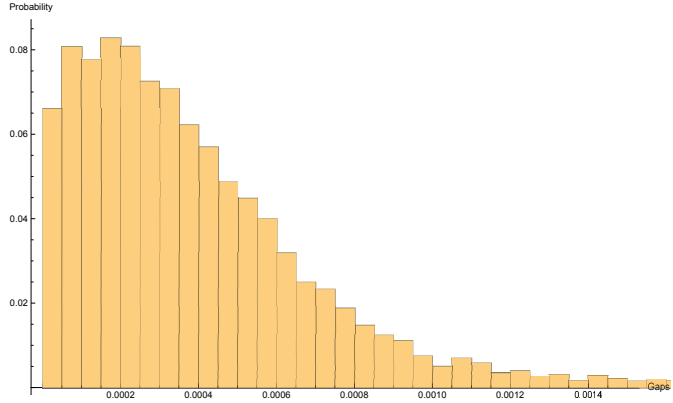}
    \caption{Eigenvalue gaps of a $20,000 \times 20,000$ matrix $\mathcal{D}_1 (A, B)$, with $A$ a random symmetric palindromic Toeplitz matrix and $B$ a random real symmetric matrix.}
    \label{fig:gaps}
\end{figure}

In \cite{MMS} it is conjectured that the limiting eigenvalue gap
distribution of symmetric palindromic Toeplitz matrices
is Poissonian, while the ensemble
of all real symmetric matrices is conjectured to have normalized
spacings given by the GOE distribution whenever the independent
matrix elements are independently chosen from a nice distribution
$p$. As $\mathcal{D}_1 (A, B)$ exhibits eigenvalue behavior
representing a hybrid of its component behaviors, we similarly
conjecture that the limiting eigenvalue gap distribution of
$\mathcal{D}_1 (A, B)$ is bounded by that of its submatrices
$A$ and $B$.

Additionally, numerical experiments in constructing $\mathcal{D}_1(A,B)$
with $A$, $B$ drawn from several pairs of random ensembles
whose limiting eigenvalue distributions are known suggests the following
conjecture; we provide support for a special case in Appendix \ref{sec:specialcaseproof}.

\begin{conj}\label{conj:LW}
Let $A$, $B$ be $N \times N$ real symmetric random matrices whose independent entries are drawn from
a fixed probability distribution $p$ with mean 0 and variance 1. Suppose
that the limiting eigenvalue distributions of $A$, $B$ have all moments
finite and appropriately bounded. Then
\begin{equation}\label{eq:conj}
    \min\left\lbrace M_k(A), M_k(B) \right\rbrace \ \leq\ M_k\left(
    \mathcal{D}_1(A,B) \right) \ \leq\ \max\left\lbrace M_k(A), M_k(B) \right\rbrace.
\end{equation}
\end{conj}

Tables \ref{tab:conj_data1} and \ref{tab:conj_data2} show experimental results
supporting Conjecture \ref{conj:LW}. We compute small moments of $\mathcal{D}_1 (A, B)$
where $A$ is either a random real symmetric matrix (i.e., a Wigner matrix)
or a random symmetric palindromic Toeplitz matrix, whose independent entries follow the standard Gaussian distribution.
The matrix $B$ is a random $3$-period block circulant matrix, whose independent entries also follow the standard Gaussian distribution (see \cite{KMMSX} for a full treatment of the construction of block-circulant matrices and their limiting eigenvalue distributions).

\begin{center}
\begin{table}[ht]
\begin{tabular}{|c|c|c|c|}\hline
Moment & $M_k(A)$ & \hspace{0.5cm} $\mathcal{D}_1(A,B)$ \hspace{0.5cm}
& $M_k(B)$ \\ \hline
4 & \ 2.000  &\  2.071    &\  2.183  \\ \hline
6 & \ 4.997  &\  5.363    &\  6.257  \\ \hline
8 & 13.985 & 15.759 & 21.974 \\ \hline
\end{tabular}
\caption{\label{tab:conj_data1} Numerical data from a $11994\times 11994$ Disco of a random real symmetric $A$ and a random $3$-period block circulant matrix $B$ supporting \eqref{conj:LW}.}
\end{table}

\begin{table}[ht]
\begin{tabular}{|c|c|c|c|}\hline
Moment & $M_k(A)$ & \hspace{0.5cm} $\mathcal{D}_1(A,B)$ \hspace{0.5cm}
& $M_k(B)$ \\ \hline
4 &\ \  2.948  &\  2.544   &\  2.330  \\ \hline
6 &\  14.863  &\  9.783    &\  7.929  \\ \hline
8 & 102.518 & 50.681 & 36.884 \\ \hline
\end{tabular}
\caption{\label{tab:conj_data2} Numerical data from a $11994\times 11994$ Disco of a random PST $A$ and a random $3$-period block circulant matrix $B$ supporting \eqref{conj:LW}.}
\end{table}
\end{center}

The computation of $M_k \left( \mathcal{D}_1(A,B) \right)$ involves the
expansion previously shown in \eqref{eq1}. The primary obstacle
is bounding the contribution of arbitrary bi-variate matrix products in
the limit as $N \to \infty$.
A central challenge
in crafting sharp bounds on the contribution of such terms is, for arbitrary
$A$ and $B$, the lack of information on the pairing configurations of
entries that do not vanish in the limit.

By the Eigenvalue Trace Lemma and \eqref{eq1}, \eqref{eq:conj} can be rewritten as
\begin{align}\label{eq:conj_1}
    \begin{split}
    & \min\left\lbrace \lim_{N\to\infty}\mathbb{E}\left[\Tr(A^k)\right], \lim_{N\to\infty}\mathbb{E}\left[\Tr(B^k)\right] \right\rbrace \\
    &\qquad \ \leq\  \ \mathbb{E}\left[\lim_{N\to\infty}\frac{1}{2^{\frac{k}{2}+1}}\Tr\left((A+B)^k+(A-B)^k\right)\right]\\
    &\qquad \ \ \ \ \ \leq\ \ \max\left\lbrace \lim_{N\to\infty}\mathbb{E}\left[\Tr(A^k)\right], \lim_{N\to\infty}\mathbb{E}\left[\Tr(B^k)\right] \right\rbrace.
    \end{split}
\end{align}

Note that \eqref{eq:conj_1} would follow immediately if, for all $N\times N$ real symmetric matrices $A$ and $B$, it were true that
\begin{align}\label{eq:conj_1_equiv}
    \min\left\lbrace \Tr(A^k), \Tr(B^k) \right\rbrace & \ \leq\ \frac{1}{2^{\frac{k}{2}+1}} \Tr\left((A+B)^k+(A-B)^k\right) \nonumber \\
    & \ \leq\ \max\left\lbrace \Tr(A^k), \Tr(B^k) \right\rbrace.
\end{align}
Unfortunately this is not the case, as evidenced by the following construction\footnote{This construction was suggested by Zhijie Chen, Jiyoung Kim, and Samuel Murray of Carnegie Mellon University. }. Let
\begin{equation}
    \mathbf{a} = \left[\begin{array}{rr}
        -33 & -31 \\
        -31 & -82
    \end{array} \right] \qquad\text{and}\qquad  \mathbf{b} = \left[\begin{array}{rr}
        26 & 78 \\
        78 & -15
    \end{array} \right]
\end{equation}
so that, for any $m \in \mathbb{Z}^+$,
\begin{equation}
    A_{2m \times 2m} = \left[ \begin{array}{ccc}
        \mathbf{a} & & \\
         & \ddots &\\
         & & \mathbf{a}
    \end{array}\right] \qquad\text{and}\qquad B_{2m \times 2m} = \left[ \begin{array}{ccc}
        \mathbf{b} & & \\
         & \ddots &\\
         & & \mathbf{b}
    \end{array}\right]
\end{equation}
are matrices of equal dimension $2m \times 2m$ with $m$ instances of $\mathbf{a}$ and $\mathbf{b}$ along their main diagonals, respectively. For $m = 10$ and $k = 4$, we compute
\begin{align}
\Tr \left(A_{20 \times 20}^4 \right) &\ =\ 889,801,750 \nonumber\\
\Tr \left(A_{20 \times 20}^4 \right) &\ =\ 869,734,090 \nonumber \\
\dfrac{\Tr\left((A_{20 \times 20} + B_{20 \times 20})^4+(A_{20 \times 20}-B_{20 \times 20})^4\right)}{2^3} & \ =\ 1,336,343,790,
\end{align}
which is clearly at odds with \eqref{eq:conj_1_equiv}.

\appendix


\section{Evidence for the Moment Conjecture}\label{sec:specialcaseproof}

We end with a proof that in the limit as the matrix sizes tend to infinity that a special case of Conjecture \ref{conj:LW} is true for the limiting moments. In case of $\D_1(A, B)$, where $A$ is PST and $B$ is Wigner, we have shown that its EESD converges in law to $(s + g)/\sqrt{2}$, where $s$ is standard free Gaussian, $g$ is standard Gaussian, and they are free. We now examine Conjecture~\ref{conj:LW} in this special case for the limiting distributions. As the limiting odd moments are zero, we only need to consider the even moments. Now \cite{BDJ} explicitly computed the even moments of the free convolution $s + g$, which shows up as the limiting spectral distribution of random Markov matrices. They showed that
\[
    m_{2k}(s + g)\ =\ \sum_{w \text{ pair-matched}} 2^{h(w)},
\]
where the height function $h(w)$ gives, for a pair-matched word $w$ of length $2k$, the number of connected pairings in $w$ (for more details on these see \cite{BDJ}). For example, $h(abab) = 0$ but $h(abba) = 2$. Clearly $h(w) \le k$ for any pair-matched word, and $h(w) = k$ for any pair-matched non-crossing word. It follows that
\[
  m_{2k}(s + g)\ \le\ 2^k \#\{\text{pair-matched words of length $2k$}\}\ =\ 2^k m_{2k}(g),
\]
and
\[
    m_{2k}(s + g)\ \ge\ \#\{\text{pair-matched non-crossing words of length $2k$}\}\ =\ 2^k m_{2k}(s).
\]
In other words,
\[
    \min\{ m_{2k}(s), m_{2k}(g) \}\ = \ m_{2k}(s) \le m_{2k}\bigg(\frac{1}{\sqrt{2}}(s + g)\bigg) \le m_{2k}(g)\ =\ \max \{ m_{2k}(s), m_{2k}(g) \}.
\]
This proves the inequalities of Conjecture~\ref{conj:LW} hold for the limits of the moments in this special case.

\bibliographystyle{alpha}

\bigskip

\end{document}